\input amstex
\documentstyle{amsppt}
\topmatter \magnification=\magstep1 \pagewidth{5.2 in}
\pageheight{6.7 in}
\abovedisplayskip=10pt \belowdisplayskip=10pt
\parskip=8pt
\parindent=5mm
\baselineskip=2pt
\title
  On the explicit formula of Euler numbers and polynomials of
  higher order
  \endtitle
\author Taekyun Kim  \endauthor

\affil{{\it Department of Mathematics Education,\\
  Kongju National University, Kongju 314-701,  Korea}\\
  {\it e-mail: {tkim$\@$kongju.ac.kr }}}\endaffil

\abstract{ In [1], the multiple Frobenius-Euler numbers and
polynomials were constructed. In this paper we give some
interesting formulae which are related to the multiple
Frobenius-Euler polynomials. The main purpose of this paper is to
give the Kummer type congruences for the multiple Frobenius-Euler
numbers.
 }\endabstract
 \thanks 2000 Mathematics Subject Classification  11S80, 11B68 \endthanks
\thanks Key words and phrases: $p$-adic $q$-integrals, Euler Numbers,
$p$-adic Euler Integral \endthanks \rightheadtext{ Frobenius-Euler
numbers and polynomials } \leftheadtext{T. Kim }
\endtopmatter

\document

\define\Z{\Bbb Z_p}

\head \S 1. Introduction \endhead Throughout this paper $\Bbb
Z_p,\,\Bbb Q_p $ and $\Bbb C_p$ will respectively  denote  the
ring of $p$-adic rational integers, the field of $p$-adic rational
numbers and the completion of algebraic closure of $\Bbb Q_p$. Let
$\nu_p$ be the normalized exponential valuation of $\Bbb C_p$ with
$|p|_p=p^{-\nu_p(p)}=p^{-1}.$ For a fixed positive integer $f$, we
set $$\align
&X=\varprojlim_N (\Bbb Z/fp^N\Bbb Z),\\
&X^{\times}=\bigcup\Sb 0<a<fp\\ (a,p)=1\endSb a+fp\Bbb Z_p,\\
&a+fp^N\Bbb Z_p=\{x\in X\mid x\equiv a\pmod{fp^N}\},
\endalign$$
where $a\in \Bbb Z$ lies in $0\leq a<fp^N, $ cf. [1, 2, 7, 8, 9,
10, 11, 12].

The classical Euler numbers were defined by
$$\frac{2}{e^t+1}=\sum_{m=0}^{\infty}\frac{E_m}{m!}t^m, \quad  (|t|<2 \pi).$$

Let $u$ be an algebraic in complex number field. Then
Frobenius-Euler polynomials were defined by
$$ \frac{1-u}{e^t-u}e^{xt}=e^{H(u,x)t}=\sum_{m=0}^{\infty}H_m(u,x)\frac{t^m}{m!}, \text{ cf. [3],}$$
where we use the technical method's notation by replacing
$H^m(u,x)$ by $H_m(u,x)$ symbolically, cf. [3, 4, 5, 6]. In the
case $x=0$, Frobenius-Euler polynomials will be called as
Frobenius-Euler numbers. We write $H_m(u)=H_m(u,0),$ cf.[3, 4, 5,
6, 7 ].
 In this paper we give some interesting
formulae which are related to the multiple Frobenius-Euler
polynomials. The main purpose of this paper is to give the Kummer
type congruences for the multiple Frobenius-Euler numbers.

\head 2. Multiple $p$-adic interpolation functions for the Euler
numbers\endhead
Let $u\neq 0$ be algebraic numbers. We fix an
embedding $\bar \Bbb Q \rightarrow \Bbb C ,$ $\bar \Bbb Q
\rightarrow \Bbb C_p,$ so that we take $u$ as an element of $\Bbb
C ,$ and $\Bbb C_p .$ The Frobenius-Euler polynomials  of order
$w$, denoted $H_n^{(w)}(u, x),$ were defined as
$$\left(\frac{1-u}{e^t-u}\right)^re^{xt}=\sum_{n=0}^{\infty}H_n^{(r)}(u,
x)\frac{t^n}{n!}, \text{ cf. [3, 4] }. \tag1$$ The values at $
x=0$ are called Frobenius-Euler numbers of order $w$; when $w=1 ,$
the polynomials or numbers are called ordinary. When $x=0$ or $w=1
,$ we often suppress that part of the notation; e.g.,
$H_m^{(r)}(u)$ denotes $H_m^{(r)}(u, 0),$ $H_n(u)$ denotes
$H_n^{(1)}(u),$ cf. [2]. By (1), we easily see that
$$\aligned
&H_n^{(w)}(u, w-x)=(-1)^nH_n^{(w)}(u^{-1}, x),\\
&H_n^{(w)}(u, x)=\sum_{i=0}^n\binom ni H_i^{(w)}(u)x^{n-i}, \text{
cf. [4]. }
\endaligned$$
where $\binom ni$ is binomial coefficient . Let $\chi$ be the
Dirichlet character of the conductor $d \in \Bbb N= \{1, 2, \cdots
\}.$ Then we define the generalized $r$-ple Euler numbers with
$\chi$ as
$$\sum_{n_1, \cdots,n_r=0}^{d-1} \frac{(1-u^d)^r\chi(\sum_{i=1}^r
x_i)\left(\prod_{i=1}^r
e^{n_it}\right)u^{rd-\sum_{i=1}^rn_i}}{(e^{dt}-u^d)^r}
=\sum_{n=0}^{\infty}H_{n, \chi}^{(r)}(u)\frac{t^n}{n!}, \text{ for
$r \in \Bbb N$} . \tag 2$$ If $\chi$ is trivial then we note that
$ H_{n, \chi}^{(r)}(u)=H_n^{(r)}(u) .$ Hence, the Eq.(2) can be
considered  as the generating function of generalized
Frobenius-Euler numbers of order $r$ in the same meaning of the
generating functions of Bernoulli numbers with order $r$, cf. [2].
By (1), and (2), we easily see that
$$H_{n, \chi}^{(r)}(u)=d^n\sum_{n_1, \cdots,
n_r=0}^{d-1}\chi(\sum_{i=1}^r n_i)u^{rd-\sum_{i=1}^rn_i}H_n^{(r)}
(u^d , \frac{n_1+\cdots+n_r}{d}) .\tag3$$ In this section we
assume that $u\in \Bbb C_p$ with $|1-u|_p \geq 1.$  For $ x\in
\Bbb Z_p,$ we say that $g$ is uniformly differentiable function at
point $ a \in \Bbb Z_p,$ and write $g \in UD(\Bbb Z_p), $ if the
difference quotients, $ F_g(x,y)=\frac{g(y)-g(x)}{y-x},$ have a
limit $ l=g^{\prime}(a)$ as $(x,y) \rightarrow (a,a)$. For $g \in
UD(\Bbb Z_p)$, let us begin with expression
$$ \frac{1}{1-u^{p^N}} \sum_{0 \leq j < p^N} u^{p^N-j}g(j)
=\sum_{0 \leq j <p^N} g(j) E_u(j+p^N \Bbb Z_p),$$ which represents
the analytic of Riemann sums for $g$ in the $p$-adic number
fields. The $p$-adic Euler integral for $g$ on $\Bbb Z_p$ is
defined as the limit of these sums (as $n \rightarrow \infty $) if
this limit exists. An invariant $p$-adic Euler integral for
function $g \in UD(\Bbb Z_p)$ is defined by
$$ \int_{\Bbb Z_p}g(x) dE_u(x)=\lim_{ N \rightarrow \infty}
\frac{1}{1-u^{p^N}}\sum_{0 \leq j < p^N} g(j) u^{p^N-j} ,\text{
 cf. [1] }. \tag 4$$
 As is well known, the generating functions of
Frobenius-Euler polynomials are represented as
$$\int_{\Bbb Z_p}e^{(x+y)t} dE_u(y)=\frac{u}{e^t-u}e^{xt}, \text{ for $t\in\Bbb Z_p$ with $|t|<p^{-\frac{1}{p-1}}$ },
\text{ cf. [1, 2] }. $$

\proclaim{Definition 1 } For $x(\neq 0), s \in \Bbb Z_p ,$ we
define Euler multiple zeta function as
$$\zeta_r(u|s, x)=\sum_{n_1, \cdots,
n_r=0}^{\infty}\frac{u^{-(n_1+\cdots+n_r)}}{(n_1+\cdots+n_r
+x)^s},\text{\rm{ cf. [1]} }. \tag5$$
\endproclaim

In [1], it is easy to see that $\zeta_r(u|-n,
x)=\frac{u^r}{(u-1)^r}H_n^{(r)}(u, x), \text{ $n \geq 0$ }. $
\proclaim{ Lemma 2} For $x=r \in \Bbb N ,$ we note that
$$\zeta_r(u|s, r)=\sum_{n_1, \cdots,
n_r=0}^{\infty}\frac{u^{-(n_1+\cdots+n_r)}}{(n_1+\cdots+n_r+r)^s}
=u^{-r}\zeta_r(u|s),$$ where $\zeta_r(u|s)=\sum_{n_1, \cdots,
n_r=1}^{\infty}\frac{u^{-(n_1+\cdots+n_r)}}{(n_1+\cdots+n_r)^s},
\text{ \rm{ cf. [1, 2]}}.$
\endproclaim

By (1), we have
$$\left(\frac{u}{u-1}\right)^rH_n^{(r)}(u, x)=p^n\sum_{n_1,
\cdots,
n_r=0}^{p-1}\frac{u^{rp-\sum_{i=1}^rn_i}}{(u^p-1)^r}H_n^{(r)}(u^p,
\frac{\sum_{i=1}^rn_i+x}{p}). \tag6$$ Applying multivariate
invariant $p$-adic Euler integral on $X ,$
$|t|_p<p^{-\frac{1}{p-1}},$ we have
$$\int_{X^r}\chi(\sum_{i=1}^r x_i)e^{t\sum_{i=1}^r x_i} dE_u(x)
=\sum_{x_1, \cdots, x_r=0}^{d-1}\frac{\chi(\sum_{i=1}^r
x_i)e^{t\sum_{i=1}^r x_i}u^{rd-\sum_{i=1}^rx_i}}{(e^{dt}-u^d)^r},
 \tag 7$$ where
$$ \int_{\Bbb X^r}f(y)dE_u(y)=\undersetbrace\text{$r$
times}\to{\int_{\Z}\int_{\Z}\cdots\int_{\Z}} f(y) dE_u(y_1)\cdots
dE_u(y_r).$$ Hence, by (7), we see that \proclaim{ Proposition 3}
For $ r\in\Bbb N ,$ we obtain the following Witt's type formula:
$$\frac{H_{n,
\chi}^{(r)}(u)}{(1-u^d)^r}=\lim_{N\rightarrow\infty}\sum_{x_1,
\cdots, x_r=0}^{dp^N-1} \frac{\chi(\sum_{i=1}^r x_i)(\sum_{i=1}^r
x_i)^n u^{rd-\sum_{i=1}^rx_i}}{(1-u^{dp^N})^r}. \tag 8$$
\endproclaim
Let $a_1, \cdots, a_r$ be nonzero $p$-adic integers. Applying
multivariate invariant $p$-adic integrals on $X,$ we have
$$\int_{X}\cdots\int_{X} e^{(\sum_{j=1}^ra_jx_j
+w)t}dE_u(x_1)\cdots dE_u(x_r)=\frac{u^r
e^{wt}}{(e^{a_1t}-u)\cdots(e^{a_rt}-u)}, \text{
$|t|_p<p^{-\frac{1}{p-1}}$}.\tag 9$$ The multiple Frobenius-Euler
polynomials in $w$ were defined by
$$\frac{(1-u)^re^{wt}}{(e^{a_1t}-u)\cdots
(e^{a_rt}-u)}=\sum_{n=0}^{\infty}H_n^{(r)}(w,u|a_1,\cdots,a_r)\frac{t^n}{n!},
\text{ cf. [3]} .\tag 10$$ Let $p$ be an odd prime number and
$\bar k=( k_1, \cdots, k_r)$ with each $k_i$ a positive integer
relative prime to $p$. The Euler-Barnes multiple zeta function was
defined as
$$\zeta_r(s,u|\alpha, \bar k)=\sum_{\mu_1,\cdots, \mu_r=0}^{\infty}\frac{u^{-(\mu_1+\cdots+\mu_r)}}
{(\alpha +k_1\mu_1+\cdots+k_r\mu_r)^s}, \text{( see [3] )}.$$ For
$n\in \Bbb N ,$ it was known that
$$\zeta_r(-n, u|\alpha, \bar k
)=\frac{u^r}{(u-1)^r}H_{n}^{(r)}(\alpha, u|k_1,\cdots, k_r).\tag11
$$ Let
$$I_0=\{ \frac{x}{p}|x=\alpha+\sum_{j=1}^rk_j i_j \equiv 0 \text{$(\mod p)$ for
some $i_1,\cdots,i_r$ with $0\leq i_1,\cdots,i_r\leq p-1$}\}. $$
By the definition of Euler-Barnes multiple zeta function, we
easily see that
$$\aligned
\zeta_r(s,u|\alpha,\bar k)&=\sum_{\Sb \mu_1,\cdots,\mu_r=0\\
\alpha+k_1\mu_1+\cdots+k_r\mu_r\equiv 0 (\mod p)\endSb}
\frac{u^{-(\mu_1+\cdots+\mu_r)}}{(\alpha+k_1\mu_1+\cdots+k_r\mu_r)^s}\\
&+
\sum_{\Sb \mu_1,\cdots,\mu_r=0\\
(\alpha+k_1\mu_1+\cdots+k_r\mu_r, p)=1\endSb}
\frac{u^{-(\mu_1+\cdots+\mu_r)}}{(\alpha+k_1\mu_1+\cdots+k_r\mu_r)^s}.\endaligned
\tag12$$ Now, we set $\mu_j=i_j+pv_j$ in thr first sum and
$\mu_j=l_j+p^{N+1}v_j$ in the second sum. Then
$$\aligned
&\zeta_r(s,u|\alpha, \bar k)=\sum_{\Sb 0\leq
i_1,\cdots,i_r<p\\\beta\in
I_0\endSb}\sum_{v_1,\cdots,v_r=0}^{\infty}\frac
{u^{-(pv_1+\cdots+pv_r+i_1+\cdots+i_r)}}
{(p\beta+p(k_1v_1+\cdots+k_rv_r))^s}\\
&+\sum_{\Sb 0\leq l_1,\cdots,l_r<p^{N+1}\\(\alpha+k_1l_1+\cdots +
k_ll_r, p)=1\endSb}\sum_{v_1,\cdots,v_r=0}^{\infty}\frac
{u^{-p^{N+1}(v_1+\cdots+v_r)-(l_1+\cdots+l_r)}}
{(\alpha+k_1l_1+\cdots +
k_ll_r+p^{N+1}(k_1v_1+\cdots+k_rv_r))^s}\\
&=\frac{1}{p^s}\sum_{\Sb 0\leq i_1,\cdots,i_r<p\\\beta\in
I_0\endSb}u^{-(i_1+\cdots+i_r)}\zeta_r(s,u^p|\beta, \bar k)\\
&+\sum_{\Sb 0\leq l_1,\cdots,l_r<p^{N+1}\\(\alpha+k_1l_1+\cdots +
k_ll_r,
p)=1\endSb}u^{-(l_1+\cdots+l_r)}\zeta_r(s,u^{p^{N+1}}|\alpha+k_1l_1+\cdots
+ k_ll_r, p^{N+1}\bar k).
\endaligned \tag 13$$
By the definition of the multiple Frobenius-Euler polynomials, we
get
$$H_n^{(r)}(\alpha, u|k_1,\cdots,k_r)=\sum_{l=0}^n\binom nl (H(u)k_1+\cdots+H(u)k_r)^l\alpha^{n-l}.\tag14$$
Let $A_l= (H(u)k_1+\cdots+H(u)k_r)^l.$ Then we easily see that
$$\aligned
A_l&=\sum_{l_1=0}^l\binom l
{l_1}H_{l_1}(u)k_r^{l_1}\cdots\sum_{l_{r-1}=0}^{l-l_1-\cdots-l_{r-2}}
\binom{l-l_1-\cdots-l_{r-2}}{l_{r-1}}\\
&\cdot H_{l_{r-1}}(u)k_2^{l_{r-1}}H_{l-l_{1}-\cdots-l_{r-1}}(u)
k_1^{l-l_1-\cdots-l_{r-1}}.\endaligned\tag 15$$ Let $m$ be a
positive integer, and let $\alpha$ and $N$ be nonnegative
integers. By (11) and (13), we note that
$$\aligned
&\frac{u^r}{(1-u)^r}H_n^{(r)}(\alpha, u|\bar
k)-p^n\frac{u^{pr}}{(1-u^p)^r}\sum_{\Sb 0\leq i_1,\cdots,i_r<p\\
\beta\in I_0 \endSb} H_n^{(r)}(\beta,u^p|\bar
k)u^{-(i_1+\cdots+i_r)}\\
&=\sum_{\Sb 0\leq l_1,\cdots,l_r<p^{N+1}\\
 (\alpha+k_1l_1+\cdots+k_rl_r,p)=1\endSb}u^{-\sum_{j=1}^rl_j}\frac{u^{rp^{N+1}}}{(1-u^{p^{N+1}})^r}
H_n^{(r)}(\alpha+\sum_{j=1}^r k_jl_j,u^{p^{N+1}}|p^{N+1}\bar k).
\endaligned\tag16$$
Now we shall look at the right hand side of identity of (16). By
(14), we note that
$$H_n^{(r)}(x,u^{p^{N+1}}|p^{N+1}\bar k)= \sum_{l=0}^n\binom nl x^{n-l}(p^{N+1})^l
(H(u)k_1+\cdots+H(u)k_r)^l. \tag17$$ Let's set
$$\sum_{\Sb 0\leq l_1,\cdots,l_r<p^{N+1}\\
 (\alpha+k_1l_1+\cdots+k_rl_r,p)=1\endSb}u^{-\sum_{j=1}^rl_j}\frac{u^{rp^{N+1}}}{(1-u^{p^{N+1}})^r}
H_n^{(r)}(\alpha+\sum_{j=1}^r k_jl_j,u^{p^{N+1}}|p^{N+1}\bar
k)=\sum_{l=0}^n T_l(n).\tag18$$ By (16), (17) and (18), we see
that
$$T_l(n)=\binom nl
(p^{N+1})^l A_l\frac{(u^{p^{N+1}})^r}{(1-u^{p^{N+1}})^r}\sum_{\Sb
0\leq
j_1,\cdots,j_r<p^{N+1}\\(\alpha+k_1j_1+\cdots+k_rj_r,p)=1\endSb}
(\alpha+k_1j_1+\cdots+k_rj_r)^{n-l}.\tag19$$ By (16) and (19), we
easily see that
$$\frac{u^r}{(1-u)^r}H_n^{(r)}(\alpha, u|\bar
k)-p^n\frac{u^{pr}}{(1-u^p)^r}\sum_{\Sb 0\leq i_1,\cdots,i_r\leq
p-1\\\beta\in I_0\endSb} H_n^{(r)}(\beta,u^p|\bar
k)u^{-(i_1+\cdots+i_r)}=\sum_{l=0}^n T_l(n).$$ Thus we note that
$$T_l(n)\equiv 0 \text{ $(\mod p^{N+1}\Bbb Z_p)$, for each $n$.}$$
Since $|1-u^{p^{N+1}}|_p \geq 1 .$ If $m\equiv n (\mod
p^{N+1}(p-1)),$ then we have
$$T_l(n)\equiv T_l(m) (\mod p^{N+1}\Bbb Z_p).$$
Therefore we obtain the following Kummer type congruence:
\proclaim{Theorem 4 } Assume that $p\geq 2r+1$. Let $\alpha$ and
$N$ be nonnegative integers and $m$ be positive integers with
$(m,p-1)=1.$ Then we have
$$\frac{u^r}{(1-u)^r}H_n^{(r)}(\alpha,u|\bar
k)-p^n\frac{u^{rp}}{(1-u^p)^r}\sum_{\Sb o\leq i_1,\cdots,i_r\leq
p-1\\\beta\in I_0\endSb}H_n^{(r)}(\beta,u^p|\bar
k)u^{-(i_1+\cdots+i_r)}\in\Bbb Z_p,$$ and if $m\equiv n(\mod
p^{N+1}(p-1))$ the congruence
$$\aligned
&\frac{u^r}{(1-u)^r}H_n^{(r)}(\alpha,u|\bar
k)-p^n\frac{u^{rp}}{(1-u^p)^r}\sum_{\Sb o\leq i_1,\cdots,i_r\leq
p-1\\\beta\in I_0\endSb}H_n^{(r)}(\beta,u^p|\bar
k)u^{-(i_1+\cdots+i_r)}\\
&\equiv \frac{u^r}{(1-u)^r}H_m^{(r)}(\alpha,u|\bar
k)-p^m\frac{u^{rp}}{(1-u^p)^r}\sum_{\Sb o\leq i_1,\cdots,i_r\leq
p-1\\\beta\in I_0\endSb}H_m^{(r)}(\beta,u^p|\bar
k)u^{-(i_1+\cdots+i_r)}, \\
&(\mod p^{N+1}).
\endaligned$$
\endproclaim

Remark. Let $$J=\{(a_1,\cdots,a_r)|0\leq a_1,\cdots,a_r\leq p-1
\text{ and } a_1+\cdots +a_r \equiv 0( \mod p)\}.$$  By using
multivariate p-adic Euler integral, we easily see that
$$\int_{\Sb Z_p^r\\x_1+\cdots+x_r\in p\Bbb Z_p \endSb} e^{t(x_1+\cdots+x_r)}dE_u(x)
=\sum_{\Sb0\leq a_1,\cdots,a_r\leq p-1\\\beta\in J\endSb}u^{-t
\beta}e^{\beta t}(\frac{1-u^p}{e^{pt}-u^p})^r, \tag 20$$  for
$|t|<p^{-\frac{1}{p-1}}$. Thus we note that
$$\int_{\Sb Z_p^r\\x_1+\cdots+x_r\in p\Bbb Z_p
\endSb}(x_1+\cdots+x_r)^ndE_u(x)=p^n\sum_{\beta\in
J}u^{-\beta}H_n^{(r)}(\beta,u^p|1,\cdots, 1).$$ For $s\in\Bbb Z_p
, $ we now define $p$-adic interpolating function as follows:
$$\zeta_{p,r}(u,s)=\int_{\Sb \Bbb Z_p^r\\x_1+\cdots+x_r \notin p\Bbb
Z_p\endSb}(x_1+\cdots+x_r)^{-s}dE_u(x).$$ For $k\geq 0,$ we have
$$\zeta_{p,r}(u,-k)=\int_{\Bbb Z_p^r}(x_1+\cdots+x_r)^n dE_u(x)-\int_{\Sb Z_p^r\\x_1+\cdots+x_r\in p\Bbb Z_p
\endSb}(x_1+\cdots+x_r)^n dE_u(x). \tag 21$$
By (20) and (21), we obtain
$$\zeta_{p,r}(u,-k)=H_n^{(r)}(u|1,\cdots,1)-p^n\sum_{\beta\in
J}u^{-\beta}H_n^{(r)}(\beta,u^p|1,\cdots,1).$$

\enddocument